\documentclass[graybox]{svmult}

\usepackage{type1cm}        
\usepackage{makeidx}         
\usepackage{graphicx}        
\usepackage{multicol}        
\usepackage[bottom]{footmisc}

\usepackage{newtxtext}       %
\usepackage[varvw]{newtxmath}       

\makeindex             

\usepackage{subcaption} 
\usepackage{booktabs} 
\usepackage{placeins}
\usepackage{tikz}
\usepackage[show]{ed}

\usepackage{multirow}
\usepackage{soul}
\usepackage{tabularx,ragged2e}
\newcolumntype{C}{>{\Centering\arraybackslash}X} 
\usepackage{caption}
\usepackage{subcaption}
\usepackage{float}
\usepackage{xcolor}
\usepackage{pdfpages}
\usepackage{amsmath}
\usepackage[linesnumbered, ruled, vlined]{algorithm2e}
\usepackage{algpseudocode}

\newcommand{\myk}{k}
\newcommand{\mystrike}{\bgroup\markoverwith{\textcolor{red}{\rule[0.5ex]{2pt}{1pt}}}\ULon}
\usepackage{ulem} 
\usepackage{csquotes}

\newcommand{\beq}{\begin{equation}}
\newcommand{\eeq}{\end{equation}}

\usepackage{color}
\usepackage{hyperref}
\definecolor{myblue}{rgb}{0,0,0.6}
\definecolor{darkgreen}{rgb}{0,0.5,0}
\hypersetup{colorlinks=true,
linkcolor=myblue,citecolor=myblue,filecolor=myblue,urlcolor=myblue}

\definecolor{escol}{rgb}{0,0.4,0}
\definecolor{sgcol}{rgb}{0,0,0.7}
\definecolor{esnewcol}{rgb}{0,0.5,0}

\newcommand{\beqs}{\begin{equation*}}
\newcommand{\eeqs}{\end{equation*}}
\newcommand{\bit}{\begin{itemize}}
\newcommand{\eit}{\end{itemize}}
\newcommand{\ben}{\begin{enumerate}}
\newcommand{\een}{\end{enumerate}}
\newcommand{\bal}{\begin{align}}
\newcommand{\eal}{\end{align}}
\newcommand{\bals}{\begin{align*}}
\newcommand{\eals}{\end{align*}}
\newcommand{\bre}{\begin{remark}}
\newcommand{\ere}{\end{remark}}
\newcommand{\bpf}{\begin{proof}}
\newcommand{\epf}{\end{proof}}
\newcommand{\ble}{\begin{lemma}}
\newcommand{\ele}{\end{lemma}}
\newcommand{\bco}{\begin{corollary}}
\newcommand{\eco}{\end{corollary}}
\newcommand{\bex}{\begin{example}}
\newcommand{\eex}{\end{example}}
\newcommand{\bth}{\begin{theorem}}
\newcommand{\enth}{\end{theorem}}

\definecolor{yxcol}{rgb}{0,0,0}

\begin{document}

\title*{Massively parallel Schwarz methods for the high frequency Helmholtz equation}
\titlerunning{RAS-PML for high frequency Helmholtz equation}
\author{Yan~Xie\orcidID{0009-0004-2861-3383}, Shihua~Gong\orcidID{0000-0003-3650-2283}, Ivan~G.~Graham\orcidID{0000-0002-5730-676X}, Euan~A.~Spence\orcidID{0000-0003-1236-4592}, Chen-Song~Zhang\orcidID{0000-0002-2213-0899}}
\authorrunning{Y.~Xie, S.~Gong, I. G. Graham, E. A. Spence, C.-S. Zhang}

\institute{Yan Xie and Chen-Song Zhang\at SKLMS, Academy of Mathematics and Systems Science, Chinese Academy of Sciences, and School of Mathematical Sciences, University of Chinese Academy of Sciences,  Beijing, 100049, China, \email{{xieyan2021,zhangcs}@lsec.cc.ac.cn}
\and Shihua Gong \at School of Science and Engineering, The Chinese University of Hong Kong, Shenzhen, Guangdong 518172, China, \email{gongshihua@cuhk.edu.cn}
\and Ivan G. Graham and Euan A. Spence\at Department of Mathematical Sciences, University of Bath, Bath, BA2 7AY, UK, \email{{masigg,eas25}@bath.ac.uk}
}

%
%
\maketitle
\abstract{We investigate the  parallel one-level overlapping Schwarz method for solving finite element discretization of high-frequency Helmholtz equations. The resulting linear systems are large, indefinite, ill-conditioned, and complex-valued. 
  We  present a practical variant of the restricted additive Schwarz method with {\color{yxcol} Perfectly Matched Layer} transmission conditions (RAS-PML), which was originally analyzed in a theoretical setting in {\tt arXiv:2404.02156},  with some numerical experiments given in {\tt arXiv:2408.16580}.  In our algorithm, the width of the overlap
  and the additional PML layer on each subdomain is allowed to decrease with $\mathcal{O}(k^{-1} \log(k))$,    as the frequency $k \rightarrow \infty$,  
  and this is observed to  ensure  good convergence  while avoiding  excessive communication.
In experiments, the proposed method achieves $\mathcal{O}(\myk^d)$ parallel scalability under Cartesian domain decomposition and exhibits $\mathcal{O}(\myk)$ iteration counts and convergence time  for $d$-dimensional Helmholtz problems ($d = 2,3$)  as
$\myk$ increases.  In this preliminary note we  restrict to experiments on  2D  problems with constant wave speed. Details, analysis and extensions to variable wavespeed and 3D will be given in {\color{yxcol} future work}.
}
\section{The Helmholtz problem}
We consider the classical Helmholtz equation given by
 \begin{equation}
    \label{eq:helmholtz0}
    \Delta u + \frac{k^2}{c^2} u = -f\quad\text{in}\quad\mathbb{R}^d,
\end{equation}
subject to the Sommerfeld radiation condition:
\begin{equation}
    \label{eq:sommerfeld}
    \frac{\partial u}{\partial r} - \mathrm{i}k u = o\left(\frac{1}{r^{(d-1)/2}}\right),\quad r = |x| \to \infty.
\end{equation}
Here $k$ is the angular frequency, $c\in C^{\infty}(\mathbb{R}^d)$ represents  the (possibly variable) wave speed,  and $f\in L^2(\mathbb{R}^d)$ is the source term. We assume that both $f$ and $1-c$ are  compactly supported within a hyper-rectangle $\Omega_{\text{int}} := \prod_{i=1}^d(a^{i}, b^{i})$. 
To keep the presentation simple, we restrict here to constant wave speed case $c= 1$.

To formulate the problem in a bounded computational domain, we restrict~\eqref{eq:helmholtz0} to the extended domain $\Omega:=\prod_{i=1}^d(a^{i}-\kappa_{\rm g}, b^{i}+\kappa_{\rm g})$, where $\kappa_{\rm g}$ denotes the thickness of the Cartesian {\color{yxcol} PML} surrounding $\Omega_{\text{int}}$. We introduce a smooth scaling function $g\in C^{\infty}(\mathbb{R})$, satisfying
\begin{equation*}
    \begin{cases}
        g(x)=g'(x) = 0, & x \leq 0 \\
        g'(x) > 0, & x > 0 \\
        g''(x) = 0, & x \in (\kappa_{\infty}, + \infty)
    \end{cases}
\end{equation*}
The scaling function $g_i(x)$ for each direction $i$ is then given by:
\begin{equation*}
    g_i(x^{i}) = \begin{cases}
        g(x^{i}-b^{i}), & x^{i} \geq b^{i} \\
        0, & x^{i} \in (a^{i}, b^{i}) \\
        -g(a^{i}-x^{i}), & x^{i} \leq a^{i}
    \end{cases}
\end{equation*}
Using these, we define the PML-modified Laplacian as:
$
    \label{eq:pml}
    \Delta_{pml} := \sum_{i=1}^d \left(\gamma_i(x^{i})^{-1}\partial_{x^{i}}\right)^2
$
where $\gamma_i:=1+\mathrm{i}g'_i$. {\color{yxcol} This formulation allows for the analytic continuation of the solution into complex coordinates, effectively absorbing outgoing waves.} In weak formulation, the truncated Helmholtz problem with PML reads: find $u\in H^1_0(\Omega)$ such that
\begin{equation}
    \label{eq:weak}
    a(u, v):= \int_{\Omega} (D\nabla u) \cdot \nabla v - (\beta\cdot\nabla u)v - k^2 u v = \int_{\Omega} f v,\quad \forall v\in H^1_0(\Omega).
\end{equation}
where $D$ is the diagonal matrix with entries $D_{ii} = \gamma_i(x^{i})^{-2}$, and $\beta$ is {the}  vector field with components $\beta_i = \gamma_i'(x^{i})/\gamma_i(x^{i})^3$. 


\section{Restricted additive Schwarz method with PML transmission conditions}
\label{sec:ras}
The work {\cite{galkowski2024convergence}} analyzed  Schwarz methods (both additive and  multiplicative) with perfectly matched {layer} transmission conditions,  and provided theoretical convergence results for these when applied to the  high-frequency Helmholtz equation. {\color{yxcol}  In this paper, we implement a practical variant of the additive method and present numerical results to demonstrate that it is scalable to $\mathcal{O}(k^d)$ pocessors. Earlier work with lower levels of parallel scaling can be found in  \cite{taus2020sweeps,leng2022trace,dai2022multidirectional}. } To keep the presentation self-contained, we
summarize here the key components of the method.

\noindent
\textbf{Cartesian covering.} 
{\color{yxcol}We cover the computational domain {$\Omega$} with $N$ overlapping subdomains $\Omega_{j,\rm int}:=\prod_{i=1}^d (a^{i}_j, b^{i}_j), j=1,\cdots,N$ obtained by extending a non-overlapping Cartesian partition in each coordinate direction. The overlap width is $\delta$.}
We then extend each interior boundary  $\partial\Omega_{j,\rm int}\not\subset \partial\Omega$ by a PML layer to obtain $\Omega_{j} = \prod_{i=1}^d (a^{i}_j-\kappa_{j,l}^i, b^{i}_j+\kappa_{j,u}^i)$, where
\begin{equation*}
    \kappa_{j,l}^i = \begin{cases}
        0, & a^{i}_j=a^{i}-\kappa_{\rm g}\\
        \kappa, & a^{i}_j>a^{i}-\kappa_{\rm g}
    \end{cases},\quad
    \kappa_{j,u}^i = \begin{cases}
        0, & b^{i}_j=b^{i}+\kappa_{\rm g}\\
        \kappa, & b^{i}_j<b^{i}+\kappa_{\rm g}
    \end{cases}
\end{equation*}
and $\kappa$ is the thickness of the Cartesian PML on subdomains. 
For each subdomain $\Omega_j$, we define a local PML problem, using the  local PML scaling function: 
\begin{equation*}
    g_{j}^i(x^{i}) = \begin{cases}
        g(x^{i}-b^{i}_j), & x^{i} \geq b^{i}_j \\
        g_{i}(x^{i}), & x^{i} \in (a^{i}_j, b^{i}_j) \\
        -g(a^{i}_j-x^{i}), & x^{i} \leq a^{i}_j
    \end{cases}
\end{equation*}
For each  $j=1, ..., N$, we let $a_j$  denote the restriction of the sesquilinear form $a$ from \eqref{eq:weak}
    to $H^1_0(\Omega_j)$. 
To combine local solutions into a global approximation, we introduce a set of non-negative partition of unity (PoU) functions $\{\chi_j\}$ on $\Omega$ based on the cover $\{\Omega_{j,\rm int}\}$, such that $\sum_{j=1}^N \chi_j \equiv 1$, and ${\rm supp}(\chi_j) \subset \Omega_{j, \rm int}$ (i.e., $\chi_j$ vanishes on the {extra} PML of $\Omega_j$) and 
$\chi_j \equiv 1$ on the non-overlapped part of $\Omega_j$, namely  $\Omega_{j,\rm novlp} := \Omega_{j,\rm int} \backslash (\bigcup_{l\neq j} \Omega_{l,\rm int})$.

\noindent
\textbf{Restricted additive Schwarz method (RAS).} 
With the above definitions, the RAS-PML method,  defined (before discretization) in \cite{galkowski2024convergence},  is given in Algorithm \ref{alg:ras}.
\begin{algorithm}[htbp!]
    \caption{\texttt{RAS\_IterSolve}}
    \label{alg:ras}
    \KwIn{Initial guess $u^{(0)}$, source term $f$, number of subdomains $N$, domain $\Omega$}
    \KwOut{Approximate solution $u$}
    Partition $\Omega$ into Cartesian subdomains $\{\Omega_j\}_{j=1}^N$\;
    \For{$n = 0, 1, 2, \dots$ \Comment{Outer loop over iterations}}{
        \For{$j = 1, \dots, N$ \Comment{Inner loop (parallel) over subdomains}}{
            Find $\mathfrak{c}_j^{(n+1)} \in H^1_0(\Omega_j)$ such that:
            \begin{equation}
            \label{eq:helm-pml-local-ras}
            a_j(\mathfrak{c}_j^{(n+1)}, v_j) = (f, v_j) - a(u^{(n)}, v_j), \quad \forall v_j \in H^1_0(\Omega_j);
            \end{equation}
        }
        Compute the global update:
        $
        \label{eq:update-ras}
        u^{(n+1)} = u^{(n)} + \sum_j \chi_j \mathfrak{c}_j^{(n+1)};
        $\\
    }
    \Return $u^{(n+1)}$\;
\end{algorithm}

{Now let $u_h$ denote the Galerkin solution of \eqref{eq:weak}}  {in a conforming finite element space} $V_h\subset H_0^1(\Omega)$. 
Then, with  $V_{h,j} := \{ v_h \vert_{\Omega_j}: v_h \in V_h\} \cap H^1_0(\Omega_j)$, the discrete version of Algorithm~\ref{alg:ras} {for computing $u_h$} is as follows.
Given the current iterate $u^{(n)}_h\in V_h$, we compute local corrections $\mathfrak c_{h,j}^{{(n+1)}} \in V_{h,j}$ 
by solving the discrete counterpart of~\eqref{eq:helm-pml-local-ras}:
\begin{equation}
\label{eq:local-ras-discrete}
a_j (\mathfrak c_{h,j}^{(n+1)}, v_{h,j} ) = ( f, \mathcal{R}^T_{h,j} v_{h,j}) - a(u^{(n)}_{h}, \mathcal{R}^T_{h,j} v_{h,j} ), \quad \forall v_{h,j}\in V_{h,j},
\end{equation}
where $\mathcal{R}_{h,j}^\top$ denotes the extension by zero from  $V_{h,j}$ to $V_h$. {The new iterate $u^{(n+1)}_{h}$ is updated as:
\begin{equation}\label{eq:update1}
u^{(n+1)}_{h} = u^{(n)}_{h} + \sum_j \tilde{\mathcal{R}}^T_{h,j} \mathfrak c_{h,j}^{(n+1)},
\end{equation}
where $\tilde{\mathcal{R}}_{h,j}^\top$ denotes the weighted extension by $\chi_j$ from  $V_{h,j}$ to $V_h$. The corresponding preconditioner is given by
\begin{equation}
\mathcal{B}^{-1} :=  \sum_j \widetilde{\mathcal{R}}_{j}^\top
    \mathcal{A}_{j}^{-1}\mathcal{R}_{j},
\end{equation}
where $\mathcal{A}_j$ is the local operator corresponding to $a_j$.
}


As shown in \cite[\S8]{galkowski2024convergence}, this is a  restricted additive Schwarz (RAS) method,
where each local subdomain problem is  equipped with {a  PML and  a Dirichlet boundary   condition}.

\noindent
\textbf{Theoretical results.} 
Suppose the PoU $\{\chi_j\}$ is $C^{\infty}$ and the Helmholtz problem is non-trapping, {\color{yxcol} which means that all rays of geometric optics escape the domain in a finite time, avoiding closed cycles or infinite
reflections} {(see \cite[\S 1.7]{galkowski2024convergence})}. {Then} results from \cite[Theorems 1.1-1.4 and 1.6]{galkowski2024convergence} establish conditions, for which, given any $M>0$ and integer $s \geq 1$, there exist constants  $\mathcal{N} \in \mathbb{N}$ and $\myk_0>0$ (both independent of $f$) such that
\begin{align}
    \label{eq:conv}
    \Vert u - u^{{(\mathcal{N})}} \Vert_{H_\myk^s(\Omega)} \leq  \myk^{-M} \Vert u - u^{(0)} \Vert_{H_\myk^1(\Omega)},\ \text{for} \ \myk>\myk_0.
\end{align}
Here the weighted Sobolev norm is defined as
$
||v||^2_{H^s_{\myk}(\Omega)} := \sum_{|\alpha|\leq s} || (\myk^{-1}\partial)^\alpha v||^2_{L^2(\Omega)}.
$
In particular, \eqref{eq:conv} implies that the fixed-point iterations converge super-algebraically fast in the number of iterations for sufficiently-large $\myk$, and that the rate of convergence improves as $\myk$ increases. 
These theoretical results apply on arbitrary {overlap $\delta$  and PML width $\kappa$,} {but these have to remain fixed as  $\myk$ increases}. 

\noindent
\textbf{Mesh refinement and number of subdomains.} 
To resolve the oscillatory solutions of {\eqref{eq:helmholtz0}}, 
the mesh size must decrease at least as {fast as}
$h = \mathcal{O}(\myk^{-1})$ as $k \rightarrow \infty$, leading to  a finite element system with at least $\mathcal{O}(\myk^d)$ degrees of freedom (DoFs). To achieve efficient parallelization,   we {partition} the domain {first} into $\mathcal{O}(\myk)$ {non-overlapping} subdomains along each coordinate direction, so that the number of DoFs within each  remains approximately constant as $\myk$ grows. {Then we add the overlap and PML layers of thickness $\delta, \kappa$ as described above}. The case of meshes refined to avoid the pollution effect {will also be discussed} in {\color{yxcol} future work}.  

Although the theory in   \cite{galkowski2024convergence}  assumes that $\delta, \kappa$ should be fixed with respect to $k$, the results here   show that these parameters can be chosen to decrease quickly as $k$ increases, leading to a scalable algorithm with no loss of convergence rate. (Preliminary  experiments were given in \cite{galkowski2024schwarz}.)   

\section{Practical improvements}
\label{sec:practical}
We propose several practical improvements to Algorithm~\ref{alg:ras} which will be illustrated by the numerical experiments below.
These include:
\begin{itemize}
\item[(i)] {\bf Reducing communication by exploiting the sparsity of the residual term.} \ It can be shown 
 that the right-hand side of \eqref{eq:local-ras-discrete} is nonzero only in overlapping regions, and thus
    the communication  of the local residuals can be restricted to these regions.
\item[(ii)] {\bf Combining PML and impedance boundary conditions for improved robustness.}\  {The fact that the local corrections  \eqref{eq:local-ras-discrete} are computed in $V_{h,j} \subset H^1_0(\Omega_j)$ implies that a Dirichlet condition is applied at the boundary $\partial \Omega_j$. (This is {a common} set-up when PML is used.) However it is simple to  apply instead the
 impedance boundary condition }
{\begin{equation}
\label{eq:impedance}
{\partial_{\nu_j}  u} - \mathrm{i}k u = 0 \quad \text{on} \quad \partial \Omega_j, 
\end{equation}}
{(as a natural boundary condition on subdomains)},  
where
$\partial_{\nu_j}$ denotes is the outward normal derivative on $\partial \Omega_j$, i.e. we use a hybrid of the PML with  impedance boundary conditions. {The local sesquilinear forms $a_j$ then  incorporate \eqref{eq:impedance} as a natural boundary condition.} This improves robustness  
when the PML is thin (see Table~\ref{tab:ras-pml-imp}).
This hybrid treatment, also used for Maxwell’s equations~\cite{collino1998perfectly}, yields smaller errors than PML {with Dirichlet boundary condition}.

\item[(iii)] {\bf Scaling of width of overlap and PML layers as $\myk$  increases.}\ 
  To balance communication and convergence, we have found that it is advantageous to let the number of grid points in the overlap and PML layers grow logarithmically with $k$. Since the mesh diameter in these experiments is of the order of a wavelength
    $\lambda= 2\pi/k$, we choose overlap $\delta$ and PML width $\kappa$ according to  the formulae 
    \begin{equation} 
      \label{igg_width}
      \delta = C_\delta \, \ell(k)\,  \frac{2\pi}{k}  \quad \text{and} \quad \kappa = C_\kappa\,  \ell(k)\,  \frac{2\pi}{k}, 
    \end{equation}
    where  $\ell(k) := \frac{k}{k-k_0} \log_2(k/k_0)$,
    {for some  constants $k_0$,  $C_\delta$ and $C_\kappa$ to be chosen}.  
  \end{itemize}


\section{Numerical experiments}
\label{sec:experiments}
\vspace{-8pt}
In our numerical tests, $\Omega$ is the unit square, and the source is {the smoothed delta function}:
$
f(x) =
    \frac{16 \myk^2}{\pi^3} \exp\left(-\frac{16 \myk^2 (x-x_c)^2}{\pi^2}\right),
    $ with the source point $x_c = (0.5, 0.5)$, {the global PML width is $\kappa_g = 3\frac{2\pi}{\myk}$ 
(i.e.,     3 wavelengths)} and the PML coordinate scaling function is taken to be $g(x) = 10 \myk x^3$.

    {For this experiment the problem \eqref{eq:helmholtz0} is discretized on a  uniform square mesh}  with 12 grid points per wavelength using bilinear elements. 
 {The domain decomposition consists of  $N$ uniform overlapping square subdomains}  and the PoU functions {are taken to be}  the {tensor}  products of the 1D PoU functions {which vary}  linearly  across the overlap regions in each direction.   
In the tables, {{\tt ovlp} and {\tt pml}}
{denote}  the number of grid points {in the   overlapping and PML regions respectively}. 
{For the formula \eqref{igg_width} of $\delta$ and $\kappa$, we use $\ell(k)$ with $k_0=150$.} In all cases, the number of processors is equal to the total number of subdomains {$N = \mathcal{O}(\myk^{2})$, chosen so that each {non-overlapping} subdomain has a bounded number of degrees of freedom as $k \rightarrow \infty$}. {Iterations are terminated when the relative residual   is below  rtol=1E-10 or the  number of iterations exceeds 500}.

\vspace{-10pt}
\subsection{Different types of boundary conditions}
\vspace{-8pt}
\label{sec:ras-pml-imp}
We {first} compare the performance of three {subdomain} boundary condition strategies:  
{RAS-PML-Drch means that a PML is combined with a  subdomain  Dirichlet boundary condition, as implied by \eqref{eq:local-ras-discrete}}.
{RAS-PML-Imp means that a PML is combined with a subdomain impedance  boundary condition, as described in \S3(ii), while  RAS-Imp means that the  PML is discarded and the
  impedance boundary condition \eqref{eq:impedance} is imposed directly on subdomain boundaries.}
{A comparison of these three strategies is given in Table~\ref{tab:ras-pml-imp}.} RAS-PML-Imp exhibits superior performance over RAS-PML-Drch, particularly at high frequencies, while RAS-Imp fails to converge {at all for high $k$}.

To further enhance convergence, we use the {\color{yxcol} Algorithm \ref{alg:ras}} as a preconditioner for the Krylov method GMRES. The results in Table~\ref{tab:ras-pml-imp-gmres} show a significant improvement 
{over Table \ref{tab:ras-pml-imp},}
with RAS-PML-Imp remaining the most effective. However, since GMRES entails higher communication costs and offers little advantage when convergence is already {satisfactory}, we only employ {the}  simple Richardson iteration
{\eqref{eq:local-ras-discrete}, \eqref{eq:update1}} in the {following} experiments. {From now on}, we refer to RAS-PML-Imp simply as RAS-PML for brevity.
\vspace{-10pt}
\begin{table}[H]
    \centering
    \caption{The performance of RAS-PML-Imp, RAS-PML-Drch and RAS-Imp.}
    \label{tab:ras-pml-imp}
    \resizebox{0.7\textwidth}{!}{
    \begin{tabular}{cccll|rl|rl|rl}
    \toprule
        \multicolumn{5}{c|}{{$\kappa$, $\delta$ in \eqref{igg_width} with $C_\kappa=3/8,\ C_\delta=1/6$}} & \multicolumn{2}{c|}{\textbf{RAS-PML-Imp}} & \multicolumn{2}{c|}{\textbf{RAS-PML-Drch}} & \multicolumn{2}{c}{\textbf{RAS-Imp}} \\ \hline
        $\myk$& grid   & $N$  & pml& ovlp & iter \ & relres & iter \ & relres       & iter \ & relres \\ 
         300  & $600^2$  & 4    &  8 & 4  & 12  \ & 1.88E-11 & 13 \ & 2.29E-11           & 30 \ & 6.44E-11 \\
         600  & $1200^2$ & 16   & 11 & 5  & 17  \ & 9.15E-11 & 18 \ & 7.22E-11           & 43 \ & 7.13E-11 \\
         1200 & $2400^2$ & 64   & 15 & 6  & 34  \ & 7.26E-11 & 35 \ & 8.83E-11           & 500 \ & 0.016 \\
         2400 & $4800^2$ & 256  & 18 & 8  & 77  \ & 8.98E-11 & 94 \ & 8.91E-11           & $\times$ \ & -diverged- \\
         4800 & $9600^2$ & 1024 & 22 & 10 & 196 \ & 9.12E-11 & 500\ & 1.17E-10           & $\times$ \ & -diverged- \\
         9600 & $19200^2$& 4096 & 26 & 12 & 500 \ & 4.43E-08 & $\times$ \ & -diverged-   & $\times$ \ & -diverged- \\
    \toprule
    \end{tabular}
    }
\end{table}

\begin{table}[H]
    \centering
    \caption{Use GMRES to improve the performance of RAS-PML-Imp, RAS-PML-Drch and RAS-Imp with slightly thin PML thickness.}
    \label{tab:ras-pml-imp-gmres}
    \resizebox{0.7\textwidth}{!}{
    \begin{tabular}{cccll|rl|rl|rl}
    \toprule
        \multicolumn{5}{c|}{{$\kappa$, $\delta$ in \eqref{igg_width} with $C_\kappa=3/8,\ C_\delta=1/6$}} & \multicolumn{2}{c|}{\textbf{RAS-PML-Imp}} & \multicolumn{2}{c|}{\textbf{RAS-PML-Drch}} & \multicolumn{2}{c}{\textbf{RAS-Imp}} \\ \hline
        $\myk$& grid   & $N$  & pml& ovlp & iter\ & relres & iter\  & relres   & iter \ & relres \\ 
         300  & $600^2$  & 4    &  8 & 4  & 11\  & 7.86E-11 & 12\   & 8.84E-11 & 29\  & 6.71E-11 \\ 
         600  & $1200^2$ & 16   & 11 & 5  & 19\  & 6.90E-11 & 19\   & 6.28E-11 & 41\ & 8.34E-11 \\ 
         1200 & $2400^2$ & 64   & 15 & 6  & 48\  & 9.30E-11 & 49\   & 9.31E-11 & 95\ & 8.20E-11 \\ 
         2400 & $4800^2$ & 256  & 18 & 8  & 88\  & 8.50E-11 & 91\   & 5.30E-11 & 310\ & 9.24E-11 \\ 
         4800 & $9600^2$ & 1024 & 22 & 10 & 154\ & 6.93E-11 & 173\  & 9.62E-11 & 500\ & 4.49E-10 \\ 
         9600 & $19200^2$& 4096 & 26 & 12 & 294\ & 8.90E-11 & 389\  & 9.92E-11 & 500\  & 9.24E-06 \\ 
    \toprule
    \end{tabular}
    }
\end{table}

\subsection{Choice of PML and overlapping width}
\label{sec:choice-pml}

We next compare the convergence behavior of RAS-PML under different choices of PML and overlap widths. 
Table~\ref{tab:ras-pml-ovlp} {shows that with  $\kappa$  chosen as a fixed multiple of wavelength, RAS-PML eventually diverges,  independently of whether $\delta$ contains a fixed or logarithmically growing number of grid-points.}
{In contrast, if $\kappa$ is chosen as in \eqref{igg_width} then convergence is obtained both for $\delta$ chosen as a multiple of wavelength or growing more quickly as in \eqref{igg_width}, with the latter producing the best iteration counts, in fact with close to linear growth as frequency increases.
Then  Table~\ref{tab:2d-time} shows that this strategy results in a total runtime that increases linearly  with respect to $\myk$.
} 
\begin{table}[H]
    \centering
    \caption{The performance of RAS-PML with different PML and overlapping width.}
    \label{tab:ras-pml-ovlp}
    \resizebox{\textwidth}{!}{
    \begin{tabular}{ccc|llr|llr|llrl|llrl}
    \toprule
        \multicolumn{3}{c|}{} & \multicolumn{6}{c|}{$\kappa=2\frac{2\pi}{\myk}$} & \multicolumn{8}{c}{{$\kappa$ in \eqref{igg_width} with $C_\kappa=1/2$}} \\ \toprule
        $\myk$ & grid &$N$ &  \multicolumn{3}{c|}{\textbf{$\delta=\frac{2}{3}\frac{2\pi}{\myk}$}} & \multicolumn{3}{c|}{{$\delta$ in \eqref{igg_width} with $C_\delta=1/6$}} & \multicolumn{4}{c|}{\textbf{$\delta=\frac{1}{3}\frac{2\pi}{\myk}$}} & \multicolumn{4}{c}{$\delta$ in \eqref{igg_width} with $C_\delta=1/6$} \\ \hline
        & & & pml & ovlp & iter & pml & ovlp & iter & pml & ovlp & iter & ratio & pml & ovlp & iter & ratio \\ \hline
         300  & $600^2$  & 4    & 25 & 8  & 6    & 25 & 4  & 7    & 12 & 4  & 7    & -    & 12 & 4  & 7    & - \\
         600  & $1200^2$ & 16   & 25 & 8  & 13   & 25 & 5  & 14   & 16 & 4  & 14   & 2.00 & 16 & 5  & 14   & 2.00 \\
         1200 & $2400^2$ & 64   & 25 & 8  & 26   & 25 & 6  & 31   & 21 & 4  & 31   & 2.21 & 21 & 6  & 31   & 2.21 \\
         2400 & $4800^2$ & 256  & 25 & 8  & 61   & 25 & 8  & 61   & 26 & 4  & 72   & 2.32 & 26 & 8  & 60   & 1.94 \\
         4800 & $9600^2$ & 1024 & 25 & 8  & 140  & 25 & 10 & 134  & 32 & 4  & 149  & 2.07 & 32 & 10 & 119  & 1.98 \\
         9600 & $19200^2$& 4096 & 25 & 8  & $\times$ & 25 & 12 & $>500$ & 38 & 4  & 298  & 2.00 & 38 & 12 & 224  & 1.88 \\
    \toprule
    \end{tabular}
    }
\end{table}
\begin{table}[H]
    \centering
    \caption{Time table where the total runtime increase with $\mathcal{O}(\myk)$.}
    \label{tab:2d-time}
    \begin{tabular}{cccll|rrr|rr}
        \toprule
        \multicolumn{5}{c|}{{$\kappa$, $\delta$ in \eqref{igg_width} with $C_\kappa=1/2,\ C_\delta=1/6$}} & \multicolumn{3}{c|}{rtol=1E-10} & \multicolumn{2}{c}{\textbf{}} \\ \hline
        $\myk$& grid   &$N$  &pml & ovlp & iter & setup(s) & solve(s) & total time(s) & total/$\myk$ \\
         300  & $600^2$  & 4    &12 & 4 & 7  & 0.95  & 0.27  & 1.22 & 0.0041\\ 
         600  & $1200^2$ & 16   &16 & 5 & 14  & 1.12  & 0.68  & 1.80 & 0.0030\\ 
         1200 & $2400^2$ & 64   &21 & 6 & 31  & 1.40  & 1.87  & 3.27 & 0.0027\\ 
         2400 & $4800^2$ & 256  &26 & 8 & 60  & 1.67  & 4.39  & 6.06 & 0.0025\\ 
         4800 & $9600^2$ & 1024 &32 & 10 & 119  & 1.94  & 11.20  & 13.15 &0.0027\\ 
         9600 & $19200^2$& 4096 &38 & 12 & 224  & 1.90  & 23.22  & 25.12 &0.0026\\ 
    \toprule
    \end{tabular}
\end{table}

\section{Conclusion}
We {describe} several improvements of our previous work {~\cite{gong2022convergence,galkowski2024convergence,galkowski2024schwarz}} to develop a practical parallel RAS-PML method for solving high-frequency Helmholtz equations.   We apply both {\color{yxcol} PMLs} and impedance boundary conditions for subdomains to make the method more robust.  {We show by experiment that allowing the PML width to contain a logarithmically growing number of grid points can yield good convergence rates  without excessive computation and communication.}   
Under a Cartesian covering with $\mathcal{O}(\myk^2)$ subdomains for 2D problems with $\mathcal{O}(\myk^2)$ DoFs, numerical experiments demonstrate that both iteration counts and total runtime grow nearly linearly for the increasing frequency $\myk$.  {Full details,  analysis and extensions to variable wavespeed and 3D are given in}~{\color{yxcol} future work}.

\noindent\textbf{Acknowledgments.} SG was supported by the National Natural Science Foundation of China (grant number 12201535) and Shenzhen Stability Science Program 2022. ES was supported by the ERC synergy grant \enquote{PSINumScat} 101167139.

\FloatBarrier
\bibliographystyle{spmpsci}

\bibliography{references}  

\end{document}